\documentclass[12pt,reqno]{amsart}
\usepackage{amsthm,amsfonts,amssymb,euscript}

\newcommand{\Z}{\ensuremath{\mathbb{Z}}}

\newcommand{\tH}{\ensuremath{\widetilde{H}}}

\numberwithin{equation}{section}
\begin{document}

\title[Complex-valued solutions of the Benjamin--Ono equation]
{Complex-valued solutions of the Benjamin--Ono equation}
\author{Alexandru D. Ionescu}
\address{University of Wisconsin--Madison}
\email{ionescu@math.wisc.edu}
\author{Carlos E. Kenig}
\address{University of Chicago}
\email{cek@math.uchicago.edu}
\thanks{The first author was
supported in part by an NSF grant, a Sloan fellowship, and a Packard fellowship.
The second author was supported in part by an NSF grant.}
\begin{abstract}
We prove that the Benjamin--Ono initial-value problem is locally well-posed for small data in the Banach spaces $\widetilde{H}^\sigma(\mathbb{R})$, $\sigma\geq 0$, of {\it{complex-valued}}  Sobolev functions with special low-frequency structure.
\end{abstract}
\maketitle

\section{Introduction}\label{section1}

In this paper we consider the Benjamin--Ono initial-value problem
\begin{equation}\label{eq-1}
\begin{cases}
\partial_tu+\mathcal{H}\partial_x^2u+\partial_x(u^2/2)=0;\\
u(0)=\phi,
\end{cases}
\end{equation}
where $\mathcal{H}$ is the Hilbert transform operator defined by
the  Fourier multiplier $-i\,\mathrm{sgn}(\xi)$. This initial-value problem has been studied extensively for real-valued data in the Sobolev spaces $H^\sigma(\mathbb{R})$, $\sigma\geq 0$ (see, for example, the introduction of \cite{IoKe} for more references). In this paper we consider small {\it{complex-valued}} data with special low-frequency structure.

We define first the Banach spaces $\widetilde{H}^\sigma(\mathbb{R})$, $\sigma\geq 0$. Let $\eta_0:\mathbb{R}\to[0,1]$ denote an even smooth function
supported in $[-8/5,8/5]$ and equal to $1$ in $[-5/4,5/4]$. For
$l\in\mathbb{Z}$ let
$\chi_l(\xi)=\eta_0(\xi/2^l)-\eta_0(\xi/2^{l-1})$, $\chi_l$
supported in $\{\xi:|\xi|\in[(5/8)\cdot2^{l},(8/5)\cdot2^{l}]\}$. For simplicity of notation, we also define the functions $\eta_l:\mathbb{R}\to[0,1]$, $l\in\mathbb{Z}$, by $\eta_l=\chi_l$ if $l\geq 1$ and $\eta_l\equiv 0$ if $l\leq -1$. We define the Banach space $B_0(\mathbb{R})$ by
\begin{equation}\label{def4'}
\begin{split}
B_0=&\{f\in L^2(\mathbb{R}):\, f \text{ supported in }[-2,2]\text{ and }\\
&||f||_{B_0}:=\inf_{f=g+h}||\mathcal{F}_{(1)}^{-1}(g)||_{L^1_x}+\sum_{k'=-\infty}^{1}2^{-k'/2}||\chi_{k'}\cdot h||_{L^2_{\xi}}<\infty\}.
\end{split}
\end{equation}
Here, and in the rest of the paper, $\mathcal{F}_{(d)}$ and $\mathcal{F}_{(d)}^{-1}$ denote the Fourier transform and the inverse Fourier transform on $\mathbb{R}^d$, $d=1,2$. For $\sigma\geq 0$ we define
\begin{equation}\label{def4}
\begin{split}
&\widetilde{H}^{\sigma}=\Big\{\phi\in L^2(\mathbb{R}):||\phi||_{\widetilde{H}^{\sigma}}^2:=||\eta_0\cdot\mathcal{F}_{(1)}(\phi)||_{B_0}^2+\sum_{k=1}^\infty 2^{2\sigma k}||\eta_k\cdot\mathcal{F}_{(1)}(\phi)||_{L^2}^2<\infty\Big\}.
\end{split}
\end{equation}
The definition shows easily that $\widetilde{H}^{\sigma}\hookrightarrow {H}^{\sigma}$, $\sigma\geq 0$, and
\begin{equation}\label{invar}
\|\delta_\lambda(\phi)\|_{\widetilde{H}^\sigma}\leq C \|\phi\|_{\widetilde{H}^\sigma}\text{ for any }\lambda\in(0,1]\text{ and }\sigma\geq 0,
\end{equation}
where $\delta_\lambda(\phi)(x):=\lambda\phi(\lambda x)$\footnote{The inequality \eqref{invar} does not improve, however, as $\lambda\to 0$, so the spaces $\widetilde{H}^\sigma$ are, in some sense, critical. Because of this we can only allow small data.}. For any Banach space $V$ and $r>0$ let $B(r,V)$ denote the open
ball $\{v\in V:||v||_V<r\}$. Let $\mathbb{Z}_+=\{0,1,\ldots\}$ and
\begin{equation*}
\widetilde{H}^\infty=\bigcap _{\sigma\in\mathbb{Z}_+}\widetilde{H}^\sigma\text{ with the induced metric}.
\end{equation*}
Our main theorem concerns local
well-posedness of
the Benjamin--Ono initial-value problem \eqref{eq-1} for small data in $\widetilde{H}^{\sigma}$, $\sigma\geq 0$.

\newtheorem{Main1}{Theorem}[section]
\begin{Main1}\label{Main1}
(a) There is a constant $\overline{\epsilon}>0$ with the property that for any $\phi\in B(\overline{\epsilon},\widetilde{H}^0)\cap \widetilde{H}^\infty$ there is a unique solution
\begin{equation*}
u=S^\infty(\phi)\in C([-1,1]:\widetilde{H}^\infty)
\end{equation*}
of the initial-value  problem \eqref{eq-1}.

(b) The mapping $\phi\to S^\infty(\phi)$ extends (uniquely) to a Lipschitz mapping
\begin{equation*}
S^0:B(\overline{\epsilon},\widetilde{H}^0)\to C([-1,1]:\widetilde{H}^0),
\end{equation*}
with the property that $S^0(\phi)$ is a solution of the initial-value problem \eqref{eq-1} for any $\phi\in B(\overline{\epsilon},\widetilde{H}^0)$ (in the sense of distributions).

(c) Moreover, for any $\sigma \in[0,\infty)$ we have the local Lipschitz bound
\begin{equation}\label{LIP}
\sup_{t\in[-1,1]}\|S^0(\phi)(t)-S^0(\phi')(t)\|_{\widetilde{H}^{\sigma}}\leq C(\sigma,R)\|\phi-\phi'\|_{\widetilde{H}^{\sigma}}
\end{equation}
for any $R>0$ and $\phi,\phi'\in B(\overline{\epsilon},\widetilde{H}^0)\cap B(R,\widetilde{H}^{\sigma})$. As a consequence, the mapping $S^0$ restricts to a locally Lipschitz mapping
\begin{equation*}
S^{\sigma}:B(\overline{\epsilon},\widetilde{H}^0)\cap \widetilde{H}^{\sigma}\to C([-1,1]:\widetilde{H}^{\sigma}).
\end{equation*}
\end{Main1}
\medskip

We discuss now some of the ingredients in the proof of Theorem \ref{Main1}. The main obstruction to proving a well-posedness result for the Benjamin--Ono equation using a fixed-point argument in some $X^{\sigma,b}$ space (in a way similar to the case of the KdV equation, see \cite{Bo}) is the lack of control of the interaction between very high and very low frequencies of solutions (cf. \cite{MoSaTz} and \cite{KoTz}). The point of the low-frequency assumption $\eta_0\cdot \widehat\phi\in B_0$ is to weaken this interaction.\footnote{Herr \cite{He} has recently used spaces similar to $\widetilde{H}^{\sigma}$ to prove local and global well-posedness for the ``dispersion-generalized'' Benjamin--Ono equation, in which the term $\mathcal{H}\partial_x^2$ is replaced by $D_x^{1+\alpha}\partial_x$, $\alpha\in(0,1]$. In this case the logarithmic divergence mentioned above does not occur.} Even with this low-frequency assumption, the use of standard $X^{\sigma,b}$ spaces for high-frequency functions (i.e. spaces defined by suitably weighted norms in the frequency space) seems to lead inevitably to logarithmic divergences (see \cite{KeCoSt} and section \ref{example}). To avoid these logarithmic divergences we work with high-frequency spaces that have two components: an $X^{\sigma,b}$-type component measured in the frequency space and a normalized $L^1_xL^2_t$ component measured in the physical space. This type of spaces have been used in the context of wave maps (see, for example, \cite{KlMa}, \cite{KlSe}, \cite{Tat1}, \cite{Tat2}, \cite{Tat3}, \cite{Tao2}, and \cite{Tao3}). Then we  prove suitable linear and bilinear estimates in these spaces and conclude the proof of Theorem \ref{Main1} using a recursive (perturbative) construction. Many of the estimates used in the proof of Theorem \ref{Main1} have already been proved in \cite{IoKe}. There are, however, several technical difficulties due to the critical definitions of the spaces $B_0$ and $\widetilde{H}^\sigma$ (see \eqref{invar}), which in this paper are larger than the corresponding spaces $B_0$ and $\widetilde{H}^\sigma$ in \cite{IoKe}. 

The rest of the paper is organized as follows: in section \ref{section2} we construct our main normed spaces and summarize some of their basic properties. In section \ref{linear} we state our main linear and bilinear estimates; most of these estimates, with the exception of Lemma \ref{Lemmac1}, are already proved in \cite{IoKe}. In section \ref{proofthm1} we combine these estimates and a recursive argument (in which we think of the nonlinear term as a perturbation) to prove Theorem \ref{Main1}. Finally, in section \ref{example} we construct two examples that justify some of the choices we make in our definitions.     

\section{The normed spaces}\label{section2}

Assume $\eta_l$ and $\chi_l$, $l\in\mathbb{Z}$, are defined as in section \ref{section1}. For $l_1\leq
l_2\in\mathbb{Z}$ let
\begin{equation*}
\eta_{[l_1,l_2]}=\sum_{l=l_1}^{l_2}\eta_l\text{ and }\eta_{\leq l_2}=\sum_{l=-\infty}^{l_2}\eta_l.
\end{equation*}
For $l\in\mathbb{Z}$ let
$I_l=\{\xi\in\mathbb{R}:|\xi|\in[2^{l-1},2^{l+1}]\}$. For
$l\in\mathbb{Z}_+$ let $\widetilde{I}_l=[-2,2]$ if
$l=0$ and $\widetilde{I}_l=I_l$ if $l\geq 1$. For $k\in\mathbb{Z}$
and $j\in\mathbb{Z}_+$ let
\begin{equation*}
\begin{cases}
&D_{k,j}=\{(\xi,\tau)\in\mathbb{R}\times\mathbb{R}:\xi\in
I_k,\,\tau-\omega(\xi)\in\widetilde{I}_j\}\,\text{ if }k\geq 1;\\
&D_{k,j}=\{(\xi,\tau)\in\mathbb{R}\times\mathbb{R}:\xi\in
I_k,\,\tau\in\widetilde{I}_j\}\,\text{ if }k\leq 0,
\end{cases}
\end{equation*}
where, for $\xi\in\mathbb{R}$,
\begin{equation}\label{omega}
\omega(\xi)=-\xi|\xi|.
\end{equation}

We define
first the normed spaces $X_k=X_k(\mathbb{R}\times\mathbb{R})$,
$k\in\mathbb{Z}_+$: for $k\geq 1$ let
\begin{equation}\label{def1}
\begin{split}
X_k=&\{f\in L^2:\, f \text{ supported in }I_k\times\mathbb{R}\text{ and }\\
&||f||_{X_k}:=\sum_{j=0}^\infty 2^{j/2}\beta_{k,j}||\eta_j(\tau-\omega(\xi))f(\xi,\tau)\,||_{L^2_{\xi,\tau}}<\infty\},
\end{split}
\end{equation}
where
\begin{equation}\label{def1'}
\beta_{k,j}=1+2^{(j-2k)/2}.
\end{equation}
The precise choice of the coefficients $\beta_{k,j}$ is important in order for all the bilinear estimates \eqref{bj1}, \eqref{bj2}, \eqref{bk1}, and \eqref{bk2} to hold (see the discussion in section \ref{example}). Notice that $2^{j/2}\beta_{k,j}\approx 2^j$ when $k$ is small. For $k=0$ we define
\begin{equation}\label{def1''}
\begin{split}
X_0=&\{f\in L^2:\, f \text{ supported in }\widetilde{I}_0\times\mathbb{R}\text{ and }\\
&||f||_{X_0}:=\sum_{j=0}^\infty\sum_{k'=-\infty}^12^{j-k'/2}||\eta_j(\tau)\chi_{k'}(\xi)f(\xi,\tau)\,||_{L^2_{\xi,\tau}}<\infty\}.
\end{split}
\end{equation}
The spaces $X_k$ are not sufficient for our purpose, due to
various logarithmic divergences involving the modulation variable. For $k\geq 100$ and $k=0$ we also define the
normed spaces $Y_k=Y_k(\mathbb{R}\times\mathbb{R})$. For $k\geq 100$ we define
\begin{equation}\label{def2}
\begin{split}
Y_k=\{f\in L^2:\,&f\text{ supported in } \bigcup_{j=0}^{k-1}D_{k,j}\text{ and }\\
&||f||_{Y_k}:=2^{-k/2}||\mathcal{F}^{-1}_{(2)}[(\tau-\omega(\xi)+i)
f(\xi,\tau)]||_{L^1_xL^2_t}<\infty\}.
\end{split}
\end{equation}
For $k=0$ we define
\begin{equation}\label{def2''}
\begin{split}
Y_0=\{f\in L^2:\,&f\text{ supported in } \widetilde{I}_0\times\mathbb{R}\text{ and }\\
&||f||_{Y_0}:=\sum_{j=0}^\infty2^j||\mathcal{F}^{-1}_{(2)}[\eta_j(\tau)
f(\xi,\tau)]||_{L^1_xL^2_t}<\infty\}.
\end{split}
\end{equation}
Then we define
\begin{equation}\label{def3}
Z_k:=X_k\text{ if }1\leq k\leq 99\text{ and }Z_k:=X_k+Y_k\text{ if }k\geq 100\text{ or }k=0.
\end{equation}
The spaces $Z_k$ are our basic normed spaces. The spaces $X_k$ are
$X^{s,b}$-type spaces; the spaces $Y_k$ are relevant due to the
local smoothing inequality
\begin{equation*}
||\partial_xu||_{L^\infty_xL^2_t}\leq
C||(\partial_t+\mathcal{H}\partial_x^2)u||_{L^1_xL^2_t}\text{ for
any }u\in\mathcal{S}(\mathbb{R}\times\mathbb{R}).
\end{equation*}

For $k\in\mathbb{Z}_+$ let
\begin{equation*}
\begin{cases}
&A_k(\xi,\tau)=\tau-\omega(\xi)+i\text{ if }k\geq 1;\\
&A_k(\xi,\tau)=\tau+i\text{ if }k=0.\\
\end{cases}
\end{equation*}
For $\sigma\geq 0$ we define the normed spaces $F^{\sigma}=F^{\sigma}(\mathbb{R}\times\mathbb{R})$ and
$N^{\sigma}=N^{\sigma}(\mathbb{R}\times\mathbb{R})$:

\begin{equation}\label{def5}
\begin{split}
F^{\sigma}&=\Big\{u\in C(\mathbb{R}:\widetilde{H}^\infty):||u||_{F^{\sigma}}^2:=\sum_{k=0}^\infty 2^{2\sigma k}
||\eta_k(\xi)(I-\partial_\tau^2)\mathcal{F}_{(2)}(u)||_{Z_k}^2
<\infty\Big\},
\end{split}
\end{equation}
and
\begin{equation}\label{def6}
\begin{split}
N^{\sigma}&=\Big\{u\in C(\mathbb{R}:\widetilde{H}^\infty):||u||_{N^{\sigma}}^2:=\sum_{k=0}^\infty 2^{2\sigma k}||\eta_k(\xi)
A_k(\xi,\tau)^{-1}\mathcal{F}_{(2)}(u)||_{Z_k}^2<\infty\Big\}.
\end{split}
\end{equation}
\medskip

We summarize now some basic properties of the spaces
$Z_k$. Using the definitions, if $k\geq 1$ and $f_k\in Z_k$ then $f_k$ can be written in the form
\begin{equation}\label{repr1}
\begin{cases}
&f_k=\sum\limits_{j=0}^\infty f_{k,j}+g_k;\\
&\sum\limits_{j=0}^\infty 2^{j/2}\beta_{k,j}||f_{k,j}||_{L^2}+||g_k||_{Y_k}\leq 2||f_k||_{Z_k},
\end{cases}
\end{equation}
such that $f_{k,j}$ is supported in $D_{k,j}$ and $g_k$ is supported in $\bigcup_{j=0}^{k-1} D_{k,j}$ (if $k\leq 99$ then $g_k\equiv 0$).
If $f_0\in Z_0$ then $f_0$ can be written in the form
\begin{equation}\label{repr2}
\begin{cases}
&f_0=\sum\limits_{j=0}^\infty \sum\limits_{k'=-\infty}^1f^{k'}_{0,j}+\sum\limits_{j=0}^\infty g_{0,j};\\
&\sum\limits_{j=0}^\infty \sum\limits_{k'=-\infty}^12^{j-k'/2}||f^{k'}_{0,j}||_{L^2}+\sum\limits_{j=0}^\infty 2^j||\mathcal{F}_{(2)}^{-1}(g_{0,j})||_{L^1_xL^2_t}\leq 2||f_0||_{Z_0},
\end{cases}
\end{equation}
such that $f^{k'}_{0,j}$ is supported in $D_{k',j}$ and $g_{0,j}$ is supported in $\widetilde{I}_0\times\widetilde{I}_j$. The  main properties of the spaces  $Z_k$ are listed  in Lemma  \ref{Lemmaa1} below (see \cite[Section 4]{IoKe} for complete proofs).

\newtheorem{Lemmaa1}{Lemma}[section]
\begin{Lemmaa1}\label{Lemmaa1}
(a) If $m,m':\mathbb{R}\to\mathbb{C}$, $k\geq 0$, and $f_k\in Z_k$ then
\begin{equation}\label{lb1}
\begin{cases}
&||m(\xi)f_k(\xi,\tau)||_{Z_k}\leq C||\mathcal{F}_{(1)}^{-1}(m)||_{L^1(\mathbb{R})}||f_k||_{Z_k};\\
&||m'(\tau)f_k(\xi,\tau)||_{Z_k}\leq C||m'||_{L^\infty(\mathbb{R})}||f_k||_{Z_k}.
\end{cases}
\end{equation}

(b) If $k\geq 1$, $j\geq 0$, and $f_k\in Z_k$ then
\begin{equation}\label{lb2}
||\eta_j(\tau-\omega(\xi))f_k(\xi,\tau)||_{X_k}\leq C||f_k||_{Z_k}.
\end{equation}

(c) If $k\geq 1$, $j\in[0,k]$, and $f_k$ is supported in $I_k\times\mathbb{R}$ then
\begin{equation}\label{lb7}
||\mathcal{F}^{-1}_{(2)}[\eta_{\leq j}(\tau-\omega(\xi))f_k(\xi,\tau)]||_{L^1_xL^2_t}\leq C||\mathcal{F}_{(2)}^{-1}(f_k)||_{L^1_xL^2_t}.
\end{equation}

(d) If $k\geq 0$, $t\in\mathbb{R}$, and $f_k\in Z_k$ then
\begin{equation}\label{lb4}
\begin{cases}
&\big|\big|\int_\mathbb{R}f_k(\xi,\tau)e^{it\tau}\,d\tau\big|\big|_{L^2_\xi}
\leq C||f_k||_{Z_k}\text{ if }k\geq 1;\\
&\big|\big|\int_\mathbb{R}f_0(\xi,\tau)e^{it\tau}\,d\tau\big|\big|_{B_0}
\leq C||f_0||_{Z_0}\text{ if }k=0.
\end{cases}
\end{equation}
As a consequence,
\begin{equation}\label{hh80}
\sup_{t\in\mathbb{R}}\|u(.,t)\|_{\widetilde{H}^\sigma}\leq C_\sigma \|u\|_{F^\sigma}\text{ for any }\sigma\geq 0\text{ and }u\in F^\sigma.
\end{equation}

(e) (maximal function estimate) If $k\geq 1$ and $(I-\partial_\tau^2)f_k\in Z_k$ then
\begin{equation}\label{pr40}
||\mathcal{F}_{(2)}^{-1}(f_k)||_{L^2_xL^\infty_t}\leq
C2^{k/2}||(I-\partial_\tau^2)f_k||_{Z_k}.
\end{equation}

(f) (local smoothing estimate) If $k\geq 1$ and $f_k\in Z_k$ then
\begin{equation}\label{lb6}
||\mathcal{F}_{(2)}^{-1}(f_k)||_{L^\infty_xL^2_t}\leq C2^{-k/2}||f_k||_{Z_k}.
\end{equation}
\end{Lemmaa1}

\section{Linear and bilinear estimates}\label{linear}

In this section we state our main linear and bilinear estimates. For any $u\in C(\mathbb{R}:L^2)$ let $\widetilde{u}(.,t)\in C(\mathbb{R}:L^2)$
denote its partial Fourier transform with respect to the variable
$x$. For $\phi\in L^2(\mathbb{R})$ let $W(t)\phi\in
C(\mathbb{R}:L^2)$ denote the solution of the free Benjamin--Ono evolution
given by
\begin{equation}\label{ni1}
[W(t)\phi]\,\,\widetilde{}\,\,(\xi,t)=e^{it\omega(\xi)}\mathcal{F}_{(1)}(\phi)(\xi),
\end{equation}
where $\omega(\xi)$ is defined in \eqref{omega}. Assume
$\psi:\mathbb{R}\to[0,1]$ is an even smooth function supported in
the interval $[-8/5,8/5]$ and equal to $1$ in the interval
$[-5/4,5/4]$. Propositions \ref{Lemmab1} and \ref{Lemmab3} below are our main linear estimates (see \cite[Section 5]{IoKe} for  proofs).

\newtheorem{Lemmab1}{Proposition}[section]
\begin{Lemmab1}\label{Lemmab1}
If $\sigma\geq 0$ and $\phi\in \widetilde{H}^{\infty}$ then
\begin{equation}\label{ak1}
||\psi(t)\cdot (W(t)\phi)||_{F^{\sigma}}\leq C_\sigma||\phi||_{\widetilde{H}^{\sigma}}.
\end{equation}
\end{Lemmab1}

\newtheorem{Lemmab3}[Lemmab1]{Proposition}
\begin{Lemmab3}\label{Lemmab3}
If $\sigma\geq 0$ and $u\in N^{\sigma}$
then
\begin{equation}\label{ak2}
\Big|\Big|\psi(t)\cdot \int_0^tW(t-s)(u(s))\,ds\Big|\Big|_{F^{\sigma}}\leq C_\sigma||u||_{N^{\sigma}}.
\end{equation}
\end{Lemmab3}
\medskip

We state now our main dyadic bilinear estimates:

\newtheorem{Lemmac1}[Lemmab1]{Lemma}
\begin{Lemmac1}\label{Lemmac1}
Assume $k\geq 20$, $k_2\in[k-2,k+2]$, $f_{k_2}\in Z_{k_2}$, and $f_{0}\in Z_{0}$. Then
\begin{equation}\label{bj1}
2^k\big|\big|\eta_k(\xi)\cdot(\tau-\omega(\xi)+i)^{-1}f_{k_2}\ast f_0\big|\big|_{Z_k}\leq C||f_{k_2}||_{Z_{k_2}}||f_{0}||_{Z_{0}}.
\end{equation}
\end{Lemmac1}

\newtheorem{Lemmac2}[Lemmab1]{Lemma}
\begin{Lemmac2}\label{Lemmac2}
Assume $k\geq 20$, $k_2\in[k-2,k+2]$, $f_{k_2}\in Z_{k_2}$, and $f_{k_1}\in Z_{k_1}$ for any $k_1\in[1,k-10]\cap\mathbb{Z}$. Then
\begin{equation}\label{bj2}
\begin{split}
2^k\big|\big|\eta_k(\xi)(\tau-&\omega(\xi)+i)^{-1}f_{k_2}\ast\sum_{k_1=1}^{k-10}f_{k_1}\big|\big|_{Z_k}\leq C||f_{k_2}||_{Z_{k_2}}\sup_{k_1\in[1,k-10]}||(I-\partial_\tau^2)f_{k_1}||_{Z_{k_1}}.
\end{split}
\end{equation}
\end{Lemmac2}

\newtheorem{Lemmae1}[Lemmab1]{Lemma}
\begin{Lemmae1}\label{Lemmae1}
Assume $k,k_1,k_2\in\mathbb{Z}_+$ have the property that
$\mathrm{max}\,(k,k_1,k_2)\leq\mathrm{min}\,(k,k_1,k_2)+30$,
$f_{k_1}\in Z_{k_1}$, and $f_{k_2}\in Z_{k_2}$. Then
\begin{equation}\label{bk1}
\big|\big|\xi\cdot \eta_k(\xi)\cdot A_k(\xi,\tau)^{-1}f_{k_1}\ast
f_{k_2}\big|\big|_{X_k}\leq
C||f_{k_1}||_{Z_{k_1}}||f_{k_2}||_{Z_{k_2}}.
\end{equation}
\end{Lemmae1}

\newtheorem{Lemmae2}[Lemmab1]{Lemma}
\begin{Lemmae2}\label{Lemmae2}
Assume $k,k_1,k_2\in\mathbb{Z}_+$, $k_1,k_2\geq k+10$, $|k_1-k_2|\leq 2$, $f_{k_1}\in Z_{k_1}$, and $f_{k_2}\in Z_{k_2}$. Then
\begin{equation}\label{bk2}
\big|\big|\xi\cdot\eta_k(\xi)\cdot A_k(\xi,\tau)^{-1}f_{k_1}\ast f_{k_2}\big|\big|_{X_k}\leq C2^{-k/4}||f_{k_1}||_{Z_{k_1}}||f_{k_2}||_{Z_{k_2}}.
\end{equation}
\end{Lemmae2}

Lemmas \ref{Lemmac2}, \ref{Lemmae1}, and  \ref{Lemmae2} are already proved in \cite[Sections 7 and 8]{IoKe} (for Lemma \ref{Lemmae1} see also the bound (8.9) in \cite{IoKe}). We only provide a proof of Lemma \ref{Lemmac1}. The main ingredient is Lemma \ref{Lemmac3} below, which follows from Lemma 7.3 in \cite{IoKe}.

\newtheorem{Lemmac3}[Lemmab1]{Lemma}
\begin{Lemmac3}\label{Lemmac3}
Assume that $k\geq 20$, $k_1\in(-\infty,1]\cap\mathbb{Z}$, $k_2\in[k-2,k+2]$, $j,j_1,j_2\in\mathbb{Z}_+$, $f_{k_1,j_1}$ is an $L^2$ function supported in $D_{k_1,j_1}$, and $f_{k_2,j_2}$ is an $L^2$ function supported in $D_{k_2,j_2}$. Then
\begin{equation}\label{bo11}
\begin{split}
2^{k}2^{j/2}\beta_{k,j}||\eta_k(\xi)&\eta_j(\tau-\omega(\xi))(\tau-\omega(\xi)+i)^{-1}(f_{k_1,j_1}\ast f_{k_2,j_2})||_{L^2}\\
&\leq C(2^{k_1/2}+2^{-k/2})^{-1}\cdot2^{j_1}||f_{k_1,j_1}||_{L^2}\cdot 2^{j_2/2}\beta_{k_2,j_2}||f_{k_2,j_2}||_{L^2}.
\end{split}
\end{equation}

If $j_1\geq k+k_1-20$ then we have the stronger bound
\begin{equation}\label{bo16}
\begin{split}
2^k&2^{j/2}\beta_{k,j}||\eta_k(\xi)\eta_j(\tau-\omega(\xi))(\tau-\omega(\xi)+i)^{-1}(f_{k_1,j_1}\ast f_{k_2,j_2})||_{L^2}\\
&\leq C2^{-\max(j,j_2)/2}(2^{k_1/2}+2^{-k/2})^{-1}\cdot 2^{j_1}||f_{k_1,j_1}||_{L^2}\cdot 2^{j_2/2}\beta_{k_2,j_2}||f_{k_2,j_2}||_{L^2}.
\end{split}
\end{equation}

In addition, $\mathbf{1}_{D_{k,j}}(\xi,\tau)(f_{k_1,j_1}\ast f_{k_2,j_2})\equiv 0$ unless
\begin{equation}\label{bo5}
\begin{cases}
&\mathrm{max}\,(j,j_1,j_2)\in[k+k_1-10,k+k_1+10]\,\text{ or}\\
&\mathrm{max}\,(j,j_1,j_2)\geq k+k_1+10\text{ and }\mathrm{max}\,(j,j_1,j_2)-\mathrm{med}\,(j,j_1,j_2)\leq 10.
\end{cases}
\end{equation}
\end{Lemmac3}

The  restriction \eqref{bo5} follows from the elementary dispersive identity
\begin{equation*}
|\omega(\xi_1+\xi_2)-\omega(\xi_1)-\omega(\xi_2)|=2\min(|\xi_1|,|\xi_2|,|\xi_1+\xi_2|)\cdot \mathrm{med}(|\xi_1|,|\xi_2|,|\xi_1+\xi_2|),
\end{equation*}
where $\mathrm{med}(\alpha_1,\alpha_2,\alpha_3)=\alpha_1+\alpha_2+\alpha_3-\max(\alpha_1,\alpha_2,\alpha_3)-\min(\alpha_1,\alpha_2,\alpha_3)$ for any $\alpha_1,\alpha_2,\alpha_3\in\mathbb{R}$.

\begin{proof}[Proof of Lemma \ref{Lemmac1}] We use the representations \eqref{repr1} and \eqref{repr2} and analyze three cases.

{\bf{Case 1:}} $f_0=f_{0,j_1}^{k_1}$ is supported in $D_{k_1,j_1}$, $f_{k_2}=f_{k_2,j_2}$ is supported in $D_{k_2,j_2}$, $j_1,j_2\geq  0$, $k_1\leq 1$, $||f_0||_{Z_0}\approx 2^{j_1-k_1/2}||f_{0,j_1}^{k_1}||_{L^2}$, and $||f_{k_2}||_{Z_{k_2}}\approx 2^{j_2/2}\beta_{k_2,j_2}||f_{k_2,j_2}||_{L^2}$. The bound \eqref{bj1} which we have to prove becomes
\begin{equation}\label{ht1}
2^k\big|\big|\eta_k(\xi)\cdot(\tau-\omega(\xi)+i)^{-1}f_{k_2,j_2}\ast f_{0,j_1}^{k_1}\big|\big|_{Z_k}\leq C2^{j_1-k_1/2}||f_{0,j_1}^{k_1}||_{L^2}\cdot2^{j_2/2}\beta_{k_2,j_2}||f_{k_2,j_2}||_{L^2}.
\end{equation}
Let $h_k(\xi,\tau)=\eta_k(\xi)(\tau-\omega(\xi)+i)^{-1}(f_{k_2,j_2}\ast f_{0,j_1}^{k_1})(\xi,\tau)$. The first observation is that for most choices of $j_1$ and $j_2$, depending on $k$ and $k_1$, the function $h_k$ is supported in a bounded number of regions $D_{k,j}$, so \eqref{bo11} suffices to control $2^k||h_k||_{X_k}$. In view of \eqref{bo5}, the function $h_k$ is supported in a bounded number of regions $D_{k,j}$, and \eqref{ht1} follows from \eqref{bo11}, unless
\begin{equation}\label{bo6}
\begin{cases}
&|j_1-(k+k_1)|\leq 10\text{ and }j_2\leq k+k_1+10\text{ or }\\
&|j_2-(k+k_1)|\leq 10\text{ and }j_1\leq k+k_1+10\text{ or }\\
&j_1,j_2\geq k+k_1-10\text{ and }|j_1-j_2|\leq 10.
\end{cases}
\end{equation}

Assume \eqref{bo6} holds. Using \eqref{bo5}, $\mathbf{1}_{D_{k,j}}(\xi,\tau)\cdot h_k\equiv 0$ unless $j\leq\max(j_1,j_2)+C$. We have two cases: if $j_1\geq k+k_1-20$, then, in view of \eqref{bo6}, $j_2\leq j_1+C$ and the function $h_k$ is supported in $\bigcup_{j\leq j_1+C}D_{k,j}$. By \eqref{bo16},
\begin{equation*}
\begin{split}
2^k&||h_k||_{X_k}\leq C2^k\sum_{j\leq j_1+C}2^{j/2}\beta_{k,j}||\eta_j(\tau-\omega(\xi))h_k(\xi,\tau)||_{L^2}\\
&\leq C\big[\sum_{j\leq j_1+C}2^{-\max(j,j_2)/2}\big]2^{-k_1/2}\cdot 2^{j_1}||f^{k_1}_{0,j_1}||_{L^2}\cdot 2^{j_2/2}\beta_{k_2,j_2}||f_{k_2,j_2}||_{L^2},
\end{split}
\end{equation*}
which suffices for \eqref{ht1}. 

Assume now that $j_1\leq k+k_1-20$, so, in view of \eqref{bo6}, $|j_2-(k+k_1)|\leq 10$ and the function $h_k$ is supported in $\bigcup_{j\leq k+k_1+C}D_{k,j}$. Then, using Lemma \ref{Lemmaa1} (b) and (c)
\begin{equation*}
\begin{split}
2^{k}||h_k||_{Z_k}&\leq C2^{k/2}||\mathcal{F}_{(2)}^{-1}[(\tau-\omega(\xi)+i)h_k(\xi,\tau)]||_{L^1_xL^2_t}\\
&\leq C2^{k/2}||\mathcal{F}_{(2)}^{-1}(f_{0,j_1}^{k_1})||_{L^2_xL^\infty_t}||\mathcal{F}_{(2)}^{-1}(f_{k_2,j_2})||_{L^2_xL^2_t}\\
&\leq C2^{(j_1-k_1)/2}||f^{k_1}_{0,j_1}||_{L^2}\cdot 2^{(k+k_1)/2}||f_{k_2,j_2}||_{L^2},
\end{split}
\end{equation*}
which suffices for \eqref{ht1} since $|j_2-(k+k_1)|\leq 10$.

{\bf{Case 2:}} $f_0=f_{0,j_1}^{k_1}$ is supported in $D_{k_1,j_1}$, $j_1\geq  0$, $k_1\leq 1$, $f_{k_2}=g_{k_2}$ is supported in $\bigcup_{j_2\leq k_2-1}D_{k_2,j_2}$, $||f_0||_{Z_0}\approx 2^{j_1-k_1/2}||f_{0,j_1}^{k_1}||_{L^2}$, and $||f_{k_2}||_{Z_{k_2}}\approx ||g_{k_2}||_{Y_{k_2}}$. The bound \eqref{bj1} which we have to prove becomes
\begin{equation}\label{ht2}
2^k\big|\big|\eta_k(\xi)\cdot(\tau-\omega(\xi)+i)^{-1}g_{k_2}\ast f_{0,j_1}^{k_1}\big|\big|_{Z_k}\leq C2^{j_1-k_1/2}||f_{0,j_1}^{k_1}||_{L^2}\cdot||g_{k_2}||_{Y_{k_2}}.
\end{equation}
We have two cases: if $j_1\geq k+k_1-20$ then let $g_{k_2,j_2}(\xi_2,\tau_2)=g_{k_2}(\xi_2,\tau_2)\eta_{j_2}(\tau_2-\omega(\xi_2))$. Using $X_k$ norms, Lemma \ref{Lemmaa1} (b), \eqref{bo5}, and \eqref{bo16}, the left-hand side of \eqref{ht2} is dominated by
\begin{equation*}
\begin{split}
&C\sum_{j,j_2\geq 0}2^{k}2^{j/2}\beta_{k,j}||\eta_k(\xi)\eta_j(\tau-\omega(\xi))(\tau-\omega(\xi)+i)^{-1}(f_{0,j_1}^{k_1}\ast g_{k_2,j_2})||_{L^2}\\
&\leq C(2^{k_1/2}+2^{-k/2})^{-1}\cdot2^{j_1}||f^{k_1}_{0,j_1}||_{L^2}\sum_{j,j_2\geq 0}2^{-\max(j,j_2)/2}\cdot 2^{j_2/2}\beta_{k_2,j_2}||g_{k_2,j_2}||_{L^2}\\
&\leq C(2^{k_1/2}+2^{-k/2})^{-1}\cdot2^{j_1}||f_{0,j_1}^{k_1}||_{L^2}\cdot ||g_{k_2}||_{Y_{k_2}},
\end{split}
\end{equation*}
which suffices to prove \eqref{ht2} in this case. 

Assume now that $j_1\leq k+k_1-20$. Let
\begin{equation*}
\begin{cases}
&g_{k_2,\mathrm{low}}(\xi_2,\tau_2)=g_{k_2}(\xi_2,\tau_2)\cdot \eta_{\leq k+k_1-20}(\tau_2-\omega(\xi_2));\\
&g_{k_2,\mathrm{high}}(\xi_2,\tau_2)=g_{k_2}(\xi_2,\tau_2)\cdot (1-\eta_{\leq k+k_1-20}(\tau_2-\omega(\xi_2))).
\end{cases}
\end{equation*}
In view of \eqref{bo5}, the function $f_{0,j_1}^{k_1}\ast g_{k_2,\mathrm{low}}$ is supported in the union of a bounded number of dyadic regions $D_{k,j}$, $|j-(k+k_1)|\leq C$. Then, using $X_k$ norms in the left-hand side of \eqref{ht2} and Lemma \ref{Lemmaa1} (c) and (f),
\begin{equation*}
\begin{split}
2^k\big|\big|\eta_k(\xi)\cdot(\tau-\omega(\xi)+i)^{-1}&g_{k_2,\mathrm{low}}\ast f_{0,j_1}^{k_1}\big|\big|_{Z_k}\leq C2^k2^{-(k+k_1)/2}||f^{k_1}_{0,j_1}\ast g_{k_2,\mathrm{low}}||_{L^2}\\
&\leq C2^{(k-k_1)/2}||\mathcal{F}_{(2)}^{-1}(f^{k_1}_{0,j_1})||_{L^2_xL^\infty_t}||\mathcal{F}_{(2)}^{-1}(g_{k_2,\mathrm{low}})||_{L^\infty_xL^2_t}\\
&\leq C2^{(k-k_1)/2}\cdot 2^{j_1/2}||f^{k_1}_{0,j_1}||_{L^2}\cdot 2^{-k/2}||g_{k_2,\mathrm{low}}||_{Y_{k_2}}\\
&\leq C2^{(j_1-k_1)/2}||f^{k_1}_{0,j_1}||_{L^2}\cdot||g_{k_2}||_{Y_{k_2}},
\end{split}
\end{equation*}
which agrees with \eqref{ht2}. To handle the part corresponding to $f_{0,j_1}^{k_1}\ast g_{k_2,\mathrm{high}}$, we notice that, in view of Lemma \ref{Lemmaa1} (b),
\begin{equation*}
\|g_{k_2,\mathrm{high}}\|_{L^2_{x,t}}\leq C2^{-(k+k_1)/2}||g_{k_2}||_{Y_{k_2}}.
\end{equation*}
Then, using $Y_k$ norms in the left-hand side of \eqref{ht2} and Lemma \ref{Lemmaa1} (c),
\begin{equation*}
\begin{split}
2^k\big|\big|\eta_k(\xi)\cdot(\tau-\omega(\xi)+i)^{-1}&g_{k_2,\mathrm{high}}\ast f_{0,j_1}^{k_1}\big|\big|_{Z_k}\leq C2^{k/2}||\mathcal{F}_{(2)}^{-1}(f^{k_1}_{0,j_1}\ast g_{k_2,\mathrm{high}})||_{L^1_xL^2_t}\\
&\leq C2^{k/2}||\mathcal{F}_{(2)}^{-1}(f^{k_1}_{0,j_1})||_{L^2_xL^\infty_t}||\mathcal{F}_{(2)}^{-1}(g_{k_2,\mathrm{high}})||_{L^2}\\
&\leq C2^{k/2}\cdot 2^{j_1/2}||f^{k_1}_{0,j_1}||_{L^2}\cdot 2^{-(k+k_1)/2}||g_{k_2}||_{Y_{k_2}}\\
&\leq C2^{(j_1-k_1)/2}||f^{k_1}_{0,j_1}||_{L^2}\cdot||g_{k_2}||_{Y_{k_2}},
\end{split}
\end{equation*}
which completes the proof of \eqref{ht2}.

{\bf{Case 3:}} $f_0=g_{0,j_1}$ is supported in $\widetilde{I}_0\times\widetilde{I}_{j_1}$, $j_1\geq 0$, $||f_0||_{Z_0}\approx 2^{j_1}||\mathcal{F}^{-1}_{(2)}(g_{0,j_1})||_{L^1_xL^2_t}$. The bound \eqref{bj1} which we have to prove becomes
\begin{equation}\label{mu1}
2^k\big|\big|\eta_k(\xi)\cdot(\tau-\omega(\xi)+i)^{-1}f_{k_2}\ast g_{0,j_1}\big|\big|_{Z_k}\leq C2^{j_1}||\mathcal{F}^{-1}(g_{0,j_1})||_{L^1_xL^2_t}\cdot||f_{k_2}||_{Z_{k_2}}.
\end{equation}
This is proved in \cite[Estimate (7.12)]{IoKe}, which completes the proof of Lemma \ref{Lemmac1}.
\end{proof}

We prove now our main bilinear estimate for functions in $F^\sigma$.

\newtheorem{Lemmat2}[Lemmab1]{Proposition}
\begin{Lemmat2}\label{Lemmat2}
If $\sigma\geq 0$ and $u,v\in F^\sigma$ then
\begin{equation}\label{vc22}
||\partial_x(uv)||_{N^\sigma}\leq C_\sigma(||u||_{F^\sigma}||v||_{F^0}+
||u||_{F^0}||v||_{F^\sigma}).
\end{equation}
\end{Lemmat2}

\begin{proof}[Proof of Proposition \ref{Lemmat2}]

For $k\in\Z_+$ we define $F_k(\xi,\tau)=\eta_k(\xi)\mathcal{F}_{(2)}(u)(\xi,\tau)$  and
$G_k(\xi,\tau)=\eta_k(\xi)\mathcal{F}_{(2)}(v)(\xi,\tau)$. Then
\begin{equation*}
\begin{cases}
&||u||^2_{F^\sigma}=\sum_{k_1=0}^\infty2^{2\sigma k_1}
||(I-\partial_\tau^2)F_{k_1}||_{Z_{k_1}}^2;\\
&||v||^2_{F^\sigma}=\sum_{k_2=0}^\infty2^{2\sigma k_2}
||(I-\partial_\tau^2)G_{k_2}||_{Z_{k_2}}^2,
\end{cases}
\end{equation*}
and
\begin{equation*}
\eta_k(\xi)\mathcal{F}[\partial_x(u\cdot v)](\xi,\tau)=
C\xi\sum_{k_1,k_2\in\mathbb{Z}_+}\eta_k(\xi)[F_{k_1}\ast
G_{k_2}](\xi,\tau).
\end{equation*}
We  observe that $\eta_k(\xi)[F_{k_1}\ast G_{k_2}](\xi,\tau)\equiv
0$ unless
\begin{equation*}
\begin{cases}
&k_1\leq k-10\text{ and }k_2\in[k-2,k+2]\,\,\text {or}\\
&k_1\in[k-2,k+2]\text{ and }k_1\leq k-10\,\,\text {or}\\
&k_1,k_2\in[k-10,k+20]\,\,\text{or}\\
&k_1,k_2\geq k+10\text{ and }|k_1-k_2|\leq 2.
\end{cases}
\end{equation*}
For $k,k_1,k_2\in\mathbb{Z}_+$ let
\begin{equation*}
H_{k,k_1,k_2}(\xi,\tau)=\eta_k(\xi)A_k(\xi,\tau)^{-1}
\xi\cdot(F_{k_1}\ast G_{k_2})(\xi,\tau).
\end{equation*}
Using the definitions,
\begin{equation}\label{nu1}
\begin{split}
||&\partial_x(u\cdot v)||_{N^{\sigma}}^2= C\sum_{k\geq 0}2^{2\sigma k}
\Big|\Big|\sum_{k_1,k_2}H_{k,k_1,k_2}\Big|\Big|_{Z_k}^2.
\end{split}
\end{equation}
For $k\in\mathbb{Z}_+$ fixed we estimate, using Lemmas
\ref{Lemmac1}, \ref{Lemmac2}, \ref{Lemmae1}, and \ref{Lemmae2},
\begin{equation*}
\begin{split}
\Big|\Big|\sum_{k_1,k_2}H_{k,k_1,k_2}&\Big|\Big|_{Z_k}\leq\sum_{|k_2-k|\leq 2}
\Big|\Big|\sum_{k_1\leq k-10}H_{k,k_1,k_2}\Big|\Big|_{Z_k}+\sum_{|k_1-k|\leq 2}
\Big|\Big|\sum_{k_2\leq k-10}H_{k,k_1,k_2}\Big|\Big|_{Z_k}\\
&+\sum_{k_1,k_2\in[k-10,k+20]}||H_{k,k_1,k_2}||_{Z_k}+
\sum_{k_1,k_2\geq k+10,\,|k_1-k_2|\leq 2}||H_{k,k_1,k_2}||_{Z_k}\\
&\leq C\big[\sum_{|k_2-k|\leq 2}||G_{k_2}||_{Z_{k_2}}\big]\cdot||u||_{F^0}
+C\big[\sum_{|k_1-k|\leq 2}||F_{k_1}||_{Z_{k_1}}\big]\cdot||v||_{F^0}\\
&+C\big[\sum_{|k_1-k|\leq 20}||F_{k_1}||_{Z_{k_1}}\big]\big[\sum_{|k_2-k|\leq 20}
||G_{k_2}||_{Z_{k_2}}\big]\\
&+C2^{-k/4}\big[\sum_{k_1\geq k}||F_{k_1}||_{Z_{k_1}}^2\big]^{1/2}
\big[\sum_{k_2\geq k}||G_{k_2}||_{Z_{k_2}}^2\big]^{1/2}.
\end{split}
\end{equation*}
The bound \eqref{vc22} follows.
\end{proof}

\section{Proof of Theorem \ref{Main1}}\label{proofthm1}

In this section we complete the proof of Theorem \ref{Main1}. The main ingredients are Propositions \ref{Lemmab1}, \ref{Lemmab3}, and \ref{Lemmat2}, and the bound \eqref{hh80}. For any interval $I=[t_0-a,t_0+a]$, $t_0\in\mathbb{R}$, $a\in[0,5/4]$, and $\sigma\geq 0$ we define the normed space
\begin{equation*}
F^\sigma(I)=\{u\in C(I:\widetilde{H}^\infty):||u||_{F^\sigma(I)}=\inf\limits_{\widetilde{u}\equiv u\text{ on }\mathbb{R}\times I}||\widetilde{u}||_{F^\sigma}<\infty\}.
\end{equation*}
With this notation, the estimate in Proposition \ref{Lemmab1} becomes
\begin{equation}\label{jm1}
||W(t-t_0)\phi||_{F^\sigma([t_0-a,t_0+a])}\leq C_\sigma ||\phi||_{\widetilde{H}^\sigma}\text{ for any }\phi\in\widetilde{H}^\infty.
\end{equation}
By combining Propositions \ref{Lemmab3} and \ref{Lemmat2} we obtain (with $I=[t_0-a,t_0+a]$)
\begin{equation}\label{jm2}
\Big|\Big|\int_{t_0}^tW(t-s)(\partial_x(u\cdot v)(s))\,ds\Big|\Big|_{F^\sigma(I)}\leq C_\sigma(||u||_{F^\sigma(I)}||v||_{F^0(I)}+||u||_{F^0(I)}||v||_{F^\sigma(I)}),
\end{equation}
for any $u,v\in F^\sigma([t_0-a,t_0+a])$, $\sigma\geq 0$. Finally, the estimate \eqref{hh80} becomes
\begin{equation}\label{jm3}
\sup_{t\in[t_0-a,t_0+a]}\|u(.,t)\|_{\widetilde{H}^\sigma}\leq C_\sigma \|u\|_{F^\sigma([t_0-a,t_0+a])}\text{ for any }u\in F^\sigma([t_0-a,t_0+a]).
\end{equation}

{\it{Proof of existence.}} We prove first the existence part of Theorem \ref{Main1}, including the Lipschitz bounds in (b) and (c). Assume, as in Theorem \ref{Main1}, that $\phi\in \widetilde{H}^\infty\cap B(\overline{\epsilon},\widetilde{H}^0)$. We define recursively $u_n\in C([-1,1]:\widetilde{H}^\infty)$, 
\begin{equation}\label{fi2}
\begin{cases}
&u_0=W(t)\phi;\\
&u_{n+1}=W(t)\phi-\frac{1}{2}\int_0^tW(t-s)(\partial_x(u_n^2)(s))\,ds\,\text{ for }\,n\in\mathbb{Z}_+.
\end{cases}
\end{equation}

We show first that
\begin{equation}\label{fi3}
\|u_n\|_{F^0([-1,1])}\leq C\|\phi\|_{\widetilde{H}^0}\text{ for any }n\in\Z_+\text{ if }\overline{\epsilon} \text{ is sufficiently small}. 
\end{equation} 
The bound \eqref{fi3} holds for $n=0$, due to \eqref{jm1}. Then, using \eqref{jm2} with $\sigma=0$, it follows that
\begin{equation*}
\|u_{n+1}\|_{F^0([-1,1])}\leq C\|\phi\|_{\widetilde{H}^0}+C\|u_n\|_{F^0([-1,1])}^2,
\end{equation*}
which leads to \eqref{fi3} by induction over $n$.

We show now that 
\begin{equation}\label{fi4}
\|u_{n}-u_{n-1}\|_{F^0([-1,1])}\leq C2^{-n}\cdot\|\phi\|_{\tH^0}\text{ for any }n\in\mathbb{Z}_+\text{ if }\overline{\epsilon} \text{ is sufficiently small}. 
\end{equation}
This is clear for $n=0$ (with $u_{-1}\equiv 0$). Then, using \eqref{jm2} with $\sigma=0$, the definition \eqref{fi2}, and \eqref{fi3}
\begin{equation*}
\|u_{n+1}-u_{n}\|_{F^0([-1,1])}\leq C\cdot \overline{\epsilon}\cdot \|u_{n}-u_{n-1}\|_{F^0([-1,1])},
\end{equation*}
which leads to \eqref{fi4} by induction over $n$.

We show now that 
\begin{equation}\label{fi5}
\|u_n\|_{F^{\sigma}([-1,1])}\leq C(\sigma, \|\phi\|_{\widetilde{H}^{\sigma}})\text{ for any }n\in\mathbb{Z}_+\text{ and }\sigma\in[0,\infty).
\end{equation}
For $\sigma\in[0,2]$, the bound \eqref{fi5} follows in the same way as the bound \eqref{fi3}, by combining \eqref{jm1}, \eqref{jm2}, and induction over $n$. Thus, for \eqref{fi5} it suffices to prove that 
\begin{equation}\label{fi5'}
\|J^{\sigma'}u_n\|_{F^{\sigma_0}([-1,1])}\leq C(\sigma', \|J^{\sigma'}\phi\|_{\widetilde{H}^{\sigma_0}})\text{ for any }n\in\mathbb{Z}_+,\,\sigma'\in\Z_+\text{ and }\sigma_0\in[0,1).
\end{equation}
We fix $\sigma_0$, and argue by induction over $\sigma'$; so we may assume that
\begin{equation}\label{fi6}
\|J^{\sigma'-1}(u_n)\|_{F^{\sigma_0}([-1,1])}\leq C(\sigma',\|J^{\sigma'-1}\phi\|_{\tH^{\sigma_0}})\text{ for any }n\in\mathbb{Z}_+,
\end{equation}
and it suffices to prove that
\begin{equation}\label{fi7}
\|\partial_{x}^{\sigma'}(u_n)\|_{F^{\sigma_0}([-1,1])}\leq C(\sigma',\|J^{\sigma'}\phi\|_{\tH^{\sigma_0}})\text{ for any }n\in\mathbb{Z}_+.
\end{equation}
The bound \eqref{fi7} for $n=0$ follows from \eqref{jm1}. We use the decomposition
\begin{equation}\label{fi8}
\partial_{x}^{\sigma'}\partial_x(u_n^2)=2\partial_x(\partial_x^{\sigma'}u_n\cdot u_n)+E_n,
\end{equation}
where
\begin{equation*}
E_n=\sum_{\sigma'_1+\sigma'_2=\sigma'\text{ and }\sigma'_1,\sigma'_2\geq 1}\partial_x(\partial_{x}^{\sigma'_1}u_n\cdot\partial_{x}^{\sigma'_2}u_n).
\end{equation*}
Using \eqref{jm2} and the induction hypothesis \eqref{fi6}, we have
\begin{equation*}
\Big|\Big|-\frac{1}{2}\int_0^tW(t-s)(E_n(s))\,ds\Big|\Big|_{F^{\sigma_0}([-1,1])}\leq C(\sigma',\|J^{\sigma'-1}\phi\|_{\tH^{\sigma_0}}).
\end{equation*}
We use now the definition \eqref{fi2}, together with the bounds \eqref{jm1} and \eqref{jm2} to conclude that
\begin{equation}\label{fi9}
\begin{split}
||&\partial_{x}^{\sigma'}u_{n+1}||_{F^{\sigma_0}([-1,1])}\leq C(\sigma',\|J^{\sigma'}\phi\|_{\tH^{\sigma_0}})\\&+C\cdot ||\partial_{x}^{\sigma'}u_n||_{F^{\sigma_0}([-1,1])}\cdot ||u_n||_{F^0([-1,1])}+C\cdot ||\partial_{x}^{\sigma'}u_n||_{F^{0}([-1,1])}\cdot ||u_n||_{F^{\sigma_0}([-1,1])}.
\end{split}
\end{equation}
Assume first that $\sigma_0=0$: the bound  \eqref{fi7} follows by induction over $n$, using  \eqref{fi3}, provided that $\overline{\epsilon}$ is sufficiently small. Thus, for any $\sigma_0\in[0,1)$, the last term in the right-hand side of \eqref{fi9} is also dominated by $C(\sigma',\|J^{\sigma'}\phi\|_{\tH^{\sigma_0}})$. Then we use again \eqref{fi9} to prove \eqref{fi7} for any $\sigma_0\in[0,1)$.

Finally, we show that
\begin{equation}\label{fi20}
\|u_{n}-u_{n-1}\|_{F^{\sigma}([-1,1])}\leq C(\sigma, \|\phi\|_{\widetilde{H}^\sigma})\cdot 2^{-n}\text{ for any }n\in\mathbb{Z}_+\text{ and }\sigma\in[0,\infty).
\end{equation}
For $\sigma\in[0,2]$, the bound \eqref{fi20} follows in the same way as the bound \eqref{fi4}, by combining \eqref{jm1}, \eqref{jm2}, induction over $n$, and the bounds \eqref{fi4} and \eqref{fi5}. Thus, for \eqref{fi20} it suffices to prove that for any $\sigma'\in\Z_+$ and $\sigma_0\in[0,1)$
\begin{equation}\label{fi20'}
\|J^{\sigma'}(u_{n}-u_{n-1})\|_{F^{\sigma_0}([-1,1])}\leq C(\sigma', \|J^{\sigma'}(\phi)\|_{\widetilde{H}^{\sigma_0}})\cdot 2^{-n}\text{ for any }n\in\mathbb{Z}_+.
\end{equation}
As before, we fix $\sigma_0$ and argue by induction over $\sigma'$. We use the decomposition
\begin{equation}\label{fi85}
\partial_x^{\sigma'}(\partial_x(u_{n}^2-u_{n-1}^2))=\partial_x[\partial_x^{\sigma'}(u_{n}-u_{n-1})\cdot (u_{n}+u_{n-1})]+\text{ other terms}.
\end{equation}
We use the induction hypothesis and \eqref{jm2} to bound the $F^{\sigma_0}([-1,1])$ norm of the other terms in the decomposition above by $C(\sigma',\|J^{\sigma'}(\phi)\|_{\widetilde{H}^{\sigma_0}})\cdot 2^{-n}$. Then, using \eqref{jm2} again, we obtain
\begin{equation}\label{fi9'}
\begin{split}
||&\partial_{x}^{\sigma'}(u_{n+1}-u_{n})||_{F^{\sigma_0}([-1,1])}\leq C(\sigma',\|J^{\sigma'}\phi\|_{\tH^{\sigma_0}})\cdot 2^{-n}\\&+C\cdot ||\partial_{x}^{\sigma'}(u_n-u_{n-1})||_{F^{\sigma_0}([-1,1])}\cdot (||u_n||_{F^0([-1,1])}+\|u_{n-1}\|_{F^0([-1,1])})\\
&+C\cdot ||\partial_{x}^{\sigma'}(u_n-u_{n-1})||_{F^{0}([-1,1])}\cdot (||u_n||_{F^{\sigma_0}([-1,1])}+||u_{n-1}||_{F^{\sigma_0}([-1,1])}).
\end{split}
\end{equation}
As before, we use this inequality first for $\sigma_0=0$, together with \eqref{fi3} and induction over $n$, to prove \eqref{fi20'} for $\sigma_0=0$. Then we incorporate the last term in the right-hand side of \eqref{fi9'} into the first term $C(\sigma',\|J^{\sigma'}\phi\|_{\tH^{\sigma_0}})\cdot 2^{-n}$, and apply the inequality again for any $\sigma_0\in[0,1)$. The bound \eqref{fi20'} follows.

We can now use \eqref{fi20} and \eqref{jm3} to construct
\begin{equation*}
u=\lim_{n\to\infty}u_n\in C([-1,1]:\tH^\infty).
\end{equation*}
In view of \eqref{fi2}, 
\begin{equation*}
u=W(t)\phi-\frac{1}{2} \int_0^tW(t-s)(\partial_x(u^2(s)))\,ds\text{ on }\mathbb{R}\times [-1,1],
\end{equation*}
so $S^\infty(\phi)=u$ is a solution of the initial-value problem \eqref{eq-1}.

For Theorem \ref{Main1} (b) and (c), it suffices to show that if $\sigma\in[0,\infty)$ and $\phi,\phi'\in B(\overline{\epsilon},\tH)\cap H^\infty$ then
\begin{equation}\label{fi80}
\sup_{t\in[-1,1]}||S^\infty(\phi)-S^\infty(\phi')||_{\tH^{\sigma}}\leq C(\sigma,||\phi||_{\tH^{\sigma}}+\|\phi'\|_{\tH^{\sigma}})\cdot ||\phi-\phi'||_{\tH^{\sigma}}.
\end{equation}
Part  (b) corresponds to the case $\sigma=0$. To prove \eqref{fi80}, we define the sequences $u_n$ and $u'_n$, $n\in\mathbb{Z}_+$, as in \eqref{fi2}. In view of \eqref{jm3}, for \eqref{fi80} it suffices to prove that for $\sigma\geq 0$ and $n\in\Z_+$
\begin{equation}\label{fi80'}
\|u_n-u'_n\|_{F^\sigma([-1,1])}\leq C(\sigma,||\phi||_{\tH^{\sigma}}+\|\phi'\|_{\tH^{\sigma}})\cdot ||\phi-\phi'||_{\tH^{\sigma}}.
\end{equation}
The proof of \eqref{fi80'} is similar to the proof of \eqref{fi20}: for $\sigma\in[0,2]$ we use \eqref{jm1}, \eqref{jm2}, and induction over $n$. For $\sigma\geq 2$ we write $\sigma=\sigma_0+\sigma'$, $\sigma'\in\Z_+$, $\sigma_0\in[0,1)$, use a decomposition similar to \eqref{fi85}, and argue by induction over $\sigma'$. This completes the proof of \eqref{LIP}.
\medskip

{\it{Proof of uniqueness.}} We prove now the uniqueness statement in Theorem \ref{Main1}: if $u_1,u_2\in C([-1,1]:\tH^\infty)$ are solutions of the initial-value problem \eqref{eq-1} and $u_1(0)=u_2(0)=\phi\in B(\overline{\epsilon},\tH^0)\cap \tH^\infty$, then $u_1\equiv u_2$. For any $T\in[0,1]$ we define
\begin{equation*}
M_l(T)=\|u_l\|_{F^0([-T,T])}\text{ for }l=1,2.
\end{equation*}
Since $u_l$, $l=1,2$, are solutions of \eqref{eq-1}, we have
\begin{equation}\label{me6}
u_l(t)=W(t)\phi-\frac{1}{2} \int_0^tW(t-s)(\partial_x(u_l^2(s)))\,ds\text{ on }\mathbb{R}\times [-T,T]\text{ for any }T\in[0,1].
\end{equation}
Thus, using \eqref{jm1} and \eqref{jm2} with $\sigma=0$,
\begin{equation*}
M_l(T)\leq C_0(\overline{\epsilon}+M_l(T)^2),\,l=1,2,\,\,T\in[0,1],
\end{equation*}
which, provided that $\overline{\epsilon}$ is sufficiently small, implies that
\begin{equation*}
M_l(T)\leq 2C_0\overline{\epsilon}\text{ or }M_l(T)\geq(2C_0)^{-1}.
\end{equation*}
The uniqueness statement would follow from \eqref{me6} and \eqref{jm2} if we could prove that
$M_l(T)\leq 2C_0\overline{\epsilon}$ for all $T\in [0,1]$, $l=1,2$. For this we need the following quasi-continuity property (see \cite[Section 12]{Tao3} for the proof of a similar statement).

\newtheorem{Lemmaw1}{Lemma}[section]
\begin{Lemmaw1}\label{Lemmaw1}
Assume $u\in C([-1,1]:\tH^\infty)$ is a solution of the initial-value problem \eqref{eq-1} and define
\begin{equation*}
M(T)=\|u\|_{F^0([-T,T])}\text{ for any }T\in[0,1]. 
\end{equation*}
Then, for any $\overline{c}>0$ there is $\overline{C}\geq 1$ (which does not depend on $u$) such that
\begin{equation}\label{me10}
\limsup_{t\to T}M(t)\leq \overline{C}\cdot \liminf_{t\to T}M(t)+\overline{c},
\end{equation}
for any $T\in[0,1]$.
\end{Lemmaw1}

\section{Two examples}\label{example}

We show first that the bilinear estimate in Lemma \ref{Lemmac2} fails logarithmically if the space $Z_k$ in the left-hand side of \eqref{bj2} is replaced with $X_k$. This is the main reason for using the spaces $Y_k$.

\newtheorem{Lemmas1}{Proposition}[section]
\begin{Lemmas1}\label{Lemmas1}
Assume $k\geq 20$. Then, for some functions $f_k\in X_k$ and $f_1\in X_1$,
\begin{equation}\label{kl1}
\begin{split}
2^k\big|\big|\eta_k(\xi)(\tau-\omega(\xi)+i)^{-1}(f_{k}\ast f_{1}) \big|\big|_{X_k}\geq  C^{-1}k||(I-\partial_\tau^2)f_{k}||_{X_{k}}||(I-\partial_\tau^2)f_{1}||_{X_{1}}.
\end{split}
\end{equation}
\end{Lemmas1}

\begin{proof}[Proof of Proposition \ref{Lemmas1}] 
With $\psi$ as in section \ref{linear}, let 
\begin{equation*}
\begin{cases}
&f_1(\xi_1,\tau_1)=\psi(10(\xi_1-2))\cdot\psi(\tau_1);\\
&f_{k}(\xi_2,\tau_2)=\psi(\xi_2-2^k)\cdot\psi(2^{-k-10}(\tau_2-\omega(\xi_2))).
\end{cases}
\end{equation*}
Then $||(I-\partial_\tau^2)f_{1}||_{X_{1}}\approx 1$ and $||(I-\partial_\tau^2)f_{k}||_{X_{k}}\approx 2^k$. An easy calculation shows that
\begin{equation*}
|(f_k\ast f_1)(\xi,\tau)|\geq C^{-1}\text{ if }\xi\in[2^{k}-1/2,2^{k}+1/2]\text{ and }|\tau-\omega(\xi)|\leq 2^k.
\end{equation*}
The bound \eqref{kl1} follows from the definitions.
\end{proof}

Our second example justifies the choice of the coefficients $\beta_{k,j}$ in \eqref{def1} and \eqref{def1'}, as well as the restriction $\sigma\geq 0$ in Proposition \ref{Lemmat2}.

\newtheorem{Lemmas2}[Lemmas1]{Proposition}
\begin{Lemmas2}\label{Lemmas2}
With the definitions of the spaces $F^\sigma$ and $N^\sigma$ in \eqref{def5} and \eqref{def6}, the inequality
\begin{equation}\label{me1}
||\partial_x(uv)||_{N^\sigma}\leq C_\sigma||u||_{F^\sigma}\cdot ||v||_{F^\sigma}
\end{equation}
holds if and only if $\sigma\geq 0$. 
\end{Lemmas2}

\begin{proof}[Proof of Proposition \ref{Lemmas2}] For $k$ large and $\psi:\mathbb{R}\to[0,1]$ as before, we define $u_+$ and $u_-$ by
\begin{equation*}
\mathcal{F}_{(2)}(u_\pm)(\xi,\tau)=\psi((\xi\mp2^k)/4)\cdot  \psi((\tau-\omega(\xi))/2^{10}).
\end{equation*}
The identity $\omega(\xi_1)+\omega(\xi-\xi_1)=-2^{k+1}\xi+O(1)$ if $|\xi|,|\xi_1-2^k|\leq C$ and an easy calculation shows thats
\begin{equation}\label{me2}
\begin{split}
\big|\eta_1(\xi)(\tau-\omega(\xi)+i)^{-1}\cdot \mathcal{F}_{(2)}[\partial_x(u_+\cdot u_-)](\xi,\tau)\big|\geq C^{-1}2^{-k}\eta_1(\xi)\cdot \psi(\tau+2^{k+1}\xi).
\end{split}
\end{equation}
Then we define $v$ by
\begin{equation*}
\mathcal{F}_{(2)}(v)(\xi,\tau)=2^{-k}\eta_1(\xi)\cdot\mathbf{1}_{[0,\infty)}(\xi)\cdot \psi((\tau+2^{k+1}\xi)/2^{10})
\end{equation*}
As before, using the identity $2^{k+1}\xi_1-\omega(\xi-\xi_1)=-\omega(\xi)+O(1)$ if $|\xi_1|,|\xi-2^k|\leq C$, we compute
\begin{equation}\label{me3}
\big|\eta_k(\xi)(\tau-\omega(\xi)+i)^{-1}\cdot \mathcal{F}_{(2)}[\partial_x(u_+\cdot v)](\xi,\tau)\big|\geq C^{-1}\psi(\xi-2^k)\cdot  \psi(\tau-\omega(\xi)).
\end{equation}
Using the definitions, we compute easily
\begin{equation*}
\|u_{\pm}\|_{F^\sigma}\approx\|u_{\pm}\|_{N^\sigma} \approx 2^{\sigma k}\text{ and }\|v\|_{F^\sigma}\approx 2^{-k/2}\cdot \beta_{1,k}.
\end{equation*}
Assuming  \eqref{me1}, we need $2^{2\sigma k}\geq C^{-1}(2^{-k/2}\beta_{1,k})$ (in view of \eqref{me2}) and 
$2^{-k/2}\beta_{1,k}\geq C^{-1}$ (in view of \eqref{me3}). This forces $\sigma\geq 0$ and $\beta_{1,k}\approx 2^{k/2}$ (compare with \eqref{def1'}.
\end{proof}

\end{document}